\documentclass[11pt,a4paper]{article}
\usepackage{amsmath,amssymb}
\newcommand{\R}{\mathbb{R}}
\newcommand{\Z}{\mathbb{Z}}
\def\cA{{\mathcal A}}
\def\cC{{\mathcal C}}

\def\cS{{\mathcal S}}
\newcommand{\ee}{\varepsilon}
\newcommand{\e}{\varepsilon}
\renewcommand{\div}{{\rm div}\,}

\newcommand{\Sup}{\displaystyle \sup}

\newcommand{\Sum}{\displaystyle \sum}

\def\bun{\bar{u}_{n} }
\def\d{\partial}
\def\ddl{\dot \Delta_l}
\def\ddj{\dot \Delta_j}
\def\ddq{\dot \Delta_q}
\def\tilde{\widetilde}
\def\hat{\widehat}

\newcommand{\D}{\Delta}

\newcommand{\n}{\nabla}
\newcommand{\N}{\frac{d}{2}}

\newcommand{\fd}{\frac{d}{2}}
\newcommand{\fdp}{\frac{d}{p}}
\newcommand{\p}{\partial}

\newcommand{\h}{\hookrightarrow}

\newcommand{\de}{\delta}

\newcommand{\dq}{\delta q}
\newcommand{\du}{\delta u}
\newcommand{\dub}{\delta \overline{u}}
\newcommand{\ub}{\overline{u}}

\newtheorem{thm}{Theorem}
\newtheorem{lem}{Lemma}

\newtheorem{prop}{Proposition}
\newtheorem{defi}{Definition}

\newtheorem{rem}{Remark}

\title{Existence of strong solutions in a larger space for the shallow-water system}

\author{Fr\'ed\'eric Charve\footnote{Universit\'e Paris-Est Cr\'eteil, Laboratoire d'Analyse et de Math\'ematiques Appliqu\'ees (UMR 8050), 61 Avenue du G\'en\'eral de Gaulle, 94 010 Cr\'eteil Cedex (France). E-mail: frederic.charve@univ-paris12.fr},  Boris Haspot  \thanks{Basque Center of Applied Mathematics, Bizkaia Technology Park, Building 500, 
E-48160, Derio (Spain), haspot@bcamath.org }}
\date{}
\begin{document}
\maketitle
\begin{abstract}
This paper is dedicated to the study of both viscous compressible barotropic fluids and Navier-Stokes equation with dependent density, when the viscosity coefficients are variable, in dimension $d\geq2$. We aim at proving the local and global well-posedness  for respectively {\it large} and \textit{small} initial data having critical Besov regularity  and  more precisely we are interested in extending the class of initial data velocity when we consider the shallow water system, improving the results in \cite{CMZ1,H2} and \cite{arma}. Our result relies on the fact that the velocity $u$ can be written as the sum of the solution $u_{L}$ of the associated linear system and a remainder velocity term $\bar{u}$; then in the specific case of the shallow-water system the remainder term $\bar{u}$ is more regular than $u_{L}$ by taking into account  the regularizing effects induced on the bilinear convection term. In particular we are able to deal with  initial velocity in $\dot{H}^{\N-1}$ as Fujita and Kato for the incompressible Navier-Stokes equations (see \cite{FK}) with an additional condition of type $u_{0}\in B^{-1}_{\infty,1}$. We would like to point out that this type of result is of particular interest when we want to deal with the problem of the convergence of the solution of compressible system to the incompressible system when the Mach number goes to $0$.
\end{abstract}
\section{Introduction}
The motion of a general barotropic compressible fluid is described by the following system:
\begin{equation}
\begin{cases}
\begin{aligned}
&\p_{t}\rho+{\rm div}(\rho u)=0,\\
&\p_{t}(\rho u)+{\rm div}(\rho u\otimes u)-{\rm div}(\mu(\rho)D(u))-\n(\lambda(\rho){\rm div} u)+\n P(\rho)=\rho f,\\
&(\rho,u)_{/t=0}=(\rho_{0},u_{0}).
\end{aligned}
\end{cases}
\label{0.1}
\end{equation}
Here $u=u(t,x)\in\R^{d}$ stands for the velocity field and $\rho=\rho(t,x)\in\R^{+}$ is the density.
The pressure $P$ is a suitable smooth function depending on the density $\rho$.
We denote by $\lambda$ and $\mu$ the two viscosity coefficients of the fluid,
which are also assumed to depend on the density and which verify some parabolic conditions for the momentum equation $\mu>0$ and $\lambda+2\mu>0$ (in the physical cases the viscosity coefficients verify $\lambda+\frac{2\mu}{d}>0$ which is a particular case of the previous assumption).
We supplement the problem with initial condition $(\rho_{0},u_{0})$ and an external force $f$.
Throughout the paper, we assume that the space variable $x\in\R^{d}$ or to the periodic
box $\mathbb{T}^{d}_{a}$ with period $a_{i}$, in the i-th direction. We restrict ourselves to the case $d\geq2$. Let us recall that in the case of constant viscosity coefficients, existence and uniqueness for (\ref{0.1}) in the case of smooth data with no vacuum has been stated in the pioneering works by Nash (see \cite{Nash}) and Matsumura and Nishida ( see \cite{MN1,MN}).\\
In this article we obtain the existence of strong solution (in finite  time with large initial data or in global time for small initial data) for an optimal class of initial data by combining two different ingredients, first the notion of scaling and seconds by taking advantage of suitable choices on  the viscosity coefficients which may confer specific structures in terms of regularity. We will detail more this last point later, but as a first example
of the importance of the viscosity coefficients we recall that in \cite{MV} A. Mellet and A. Vasseur have obtained 
the stability of the global weak solutions for the Saint-Venant system by using new entropy giving additional regularity on the gradient of the density and on the integrability of the velocity.\\
Let us recall the fundamental notion of scaling for system (\ref{0.1}). Indeed guided in our approach by numerous works dedicated to the incompressible Navier-Stokes equation (see in particular the pioneering works of Fujita and Kato concerning the existence of strong solutions for the incompressible Navier-Stokes equations in \cite{FK}), we aim at solving (\ref{0.1}) in the case where the data $(\rho_{0},u_{0})$ have \textit{critical} regularity. By critical, we mean that we want to solve the system (\ref{0.1}) in functional spaces with invariant norm by the natural changes of scales which leave (\ref{0.1}) invariant. More precisely in 
 the case of barotropic fluids,  the following transformations:
\begin{equation}
(\rho(t,x),u(t,x))\longrightarrow (\rho(l^{2}t,lx),lu(l^{2}t,lx)),\;\;\;l\in\R,
\label{1}
\end{equation}
verify this property, provided that the pressure term has been changed accordingly.
This notion of critical functional frameworks has been extensively used in order to obtain optimal class of initial data for the existence of global strong solution (see \cite{CD,CMZ,DL,DW, H1, H2} in the case of constant viscosity coefficients). In particular the first result on the existence of strong solutions in spaces invariant for the scaling of the equations (when the viscosity coefficients are constant) is due to R. Danchin in \cite{DL} when the
initial data $(q_{0}=\rho_{0}-\bar{\rho},u_{0})$ (with $\bar{\rho}>0$) are in $\dot{B}^{\frac{d}{2}}_{2,1}\times (\dot{B}^{\frac{d}{2}-1}_{2,1})^{d}$.  In \cite{DW}, R. Danchin generalizes the previous results by working with more general Besov space of the type $\dot{B}^{\fdp}_{p,1}\times (\dot{B}^{\fdp-1}_{p,1})^{d}$ with some restrictions on the choice of $p$ ($p\leq d$) due to some limitation in the application of the paraproduct law. The fact that  $p=p_{1}$ is a consequence of the  strong coupling between the density and the velocity equations, indeed the pressure term is considered as a remainder for the parabolic operator in the momentum equation of (\ref{1}). In \cite{H2}, the second author generalizes the results of \cite{DW} inasmuch as we have no restriction on the size of $p$ for the initial density. To do this we are working with a new variable, the effective velocity which allows  to cancel out the coupling between the pressure and the velocity.\\
In the case of global strong solution in critical space for small initial data, we would like to recall the works of R. Danchin in \cite{DG} who shows for the first time a result of global existence of  strong solution close to a stable equilibrium when the initial data verify $(q_{0},u_{0})\in (\dot{B}^\fd_{2,1}\cap \dot{B}^{\fd-1}_{2,1})\times \dot{B}^{\fd-1}_{2,1}$. The main difficulty is to get estimates on the linearized system given that the velocity and the density are coupled via the pressure. What is crucial in this work is the smoothing effect on the velocity and a $L^{1}$ decay on $\rho-\bar{\rho}$ (this plays a key role to control the pressure term). This work was generalized in the framework of Besov space with a Lebesgue index different of $p=2$ by Q. Chen et al in \cite{CMZ} , F. Charve and R. Danchin in \cite{CD} and B. Haspot in \cite{arma}.\\
However very few articles really take into account the structure of the viscosity coefficients, indeed most of them deal with constant viscosity coefficients or consider the system (\ref{0.1}) as a perturbation of the previous case. In addition to have a norm invariant by (\ref{1}), appropriate functional spaces for solving (\ref{0.1}) must provide a control on the $L^{\infty}$ norm of the density for at least two reasons; the first corresponds to avoid the vacuum and  to ensure the parabolicity of system (\ref{0.1}), the second is to assume some property of multiplier on the density in order to deal with the term $\frac{1}{\rho}\D $. It explains why the authors restrict their study to the case where the initial density is assumed in $\dot{B}^{\fdp}_{p,1}$ with $p\in [1,+\infty[$ suitably chosen, indeed  $\dot{B}^{\fdp}_{p,1}$  is embedded in $L^{\infty}$. Furthermore in order to propagate this regularity on the density via the mass equation, it appears necessary to assume that the velocity is Lipschitz, it means $\n u\in L^{1}_{T}(\dot{B}^{\frac{d}{p}}_{p,1})\h L^{1}_{T}(L^{\infty}) $ with $p\in [1,+\infty[$ and $T>0$. We would like to point out that this necessary Lipschitz control on the velocity seems to prevent any hope of existence of strong solutions when $u_{0}$ is only assumed to belong in $\dot{H}^{\N-1}$.\\
 Indeed in \cite{CD,CMZ,DL,DW, H1, H2}  the idea is  to propagate a  $\widetilde{L}^{\infty}(\dot{B}^{\N-1}_{2,1})\cap L^{1}(\dot{B}^{\N+1}_{2,1})$ regularity on the velocity $u$ via the regularizing effects induced by the momentum equation written in its eulerian form:
\begin{equation}
\p_{t}u-\frac{\mu}{\rho}\D u-\frac{\mu+\lambda}{\rho}\n{\rm div}u=R,
\label{parabol1}
\end{equation} 
with: $R=-u\cdot\n u-\frac{\n P(\rho)}{\rho}+f$. However if we only were interested in obtaining a regularity in $\widetilde{L}^{\infty}(\dot{H}^{\N-1})\cap L^{1}(\dot{H}^{\N+1})$ , it would be necessary to deal with the term $\frac{1}{\rho}\D u$ and getting enough regularity in order that $\frac{1}{\rho}$ remains in a multiplier  space ${\cal M}(\dot{H}^{\N-1})$ of $\dot{H}^{\N-1}$. Typically $\dot{H}^{\N}\cap L^{\infty}$ is embedded in ${\cal M}(\dot{H}^{\N-1})$. However in this case how to propagate the regularity $\dot{H}^{\N}\cap L^{\infty}$ on the density as our velocity is only assumed in $ \widetilde{L}^{1}(\dot{H}^{\N+1})$ (it means not necessary Lipshitz)?\\
We want to partially solve this question in the case of specific viscosity coefficients. As we explained above, one of the main difficulty is linked to the treatment of heat equation with variable coefficients. We would like to work with the shallow water system (i.e $\mu(\rho)=\mu\rho$ and $\lambda(\rho)=0$, and to simplify we will take $\mu=1$) that we can rewrite in the following form:
\begin{equation}
\begin{cases}
\begin{aligned}
&\p_{t}\rho+{\rm div}(\rho u)=0,\\
&\p_{t}u+u\cdot\n u-\cA u-2D(u).\n(\ln\rho)+\n(G(\rho))=0,\\
&(\rho,u)_{/t=0}=(\rho_{0},u_{0}),
\end{aligned}
\end{cases}
\end{equation}
where the operator $\mathcal{A}$ is defined by: $\mathcal{A} u =\Delta u +\n \div u$ and $G^{'}(\rho)=\frac{P^{'}(\rho)}{\rho}$. Roughly speaking, we are led to consider a basic heat equation with remainder terms, and in particular it seems possible to propagate in this case the regularity $\widetilde{L}^{\infty}(\dot{H}^{\N-1})\cap \widetilde{L}^{1}(\dot{H}^{\N+1})$. However due to the strict coupling between the pressure and the velocity, it will be necessary to control the density in $L^{1}_{T}(\dot{B}^{\N}_{2,1})$ (the third index $1$ is due to multiplier space reason). More precisely we will split the solution $u$ into the following sum $u=u_{L}+\bar{u}$ with:
\begin{equation}
\begin{cases}
\begin{aligned}
&\p_{t}u_{L}-\cA u_{L}=0\\
&u_{/t=0}=u_{0}
\end{aligned}
\end{cases}
\end{equation}
In the present article, if we assume that $u_{0}$ belongs to $H^{\N-1}\cap B^{-1}_{\infty,1}$ then, with usual methods we show that $\bar{u}$ is more regular than $u_L$ ($\bar{u}$ will be in $L^{1}_{T}(\dot{B}^{\N+1}_{2,1})$) which will help us propagating the regularity on $u$.\\
To simplify the notation, we assume from now on that $\bar{\rho}=1$. Hence as long as $\rho$ does not vanish, the equations for ($q=\rho-1$,$u$) read:
\begin{equation}
\begin{cases}
\begin{aligned}
&\p_{t}q+u\cdot\n q+(1+q){\rm div}u=0,\\
&\p_{t}u+u\cdot\n u-{\cal A}u-2D(u).\n(\ln(1+q))+\n(G(1+q))=0,
\end{aligned}
\end{cases}
\label{0.6}
\end{equation}
We can now state our main result:
\begin{thm}\sl{
Let $P$ be a suitably smooth function of the density. Assume that $u_{0}\in \dot{B}^{\fd-1}_{2,2}\cap \dot{B}^{-1}_{\infty,1}$, $\div u_{0}\in \dot{B}^{\fd-2}_{2,1}$, $q_{0}\in \dot{B}^{\fd-1}_{2,1}$ and that there exists $c>0$ such that $\rho_{0}\geq c$. Then there exists a time $T$ such that there exists a unique solution $(q,u)$ for system (\ref{0.1}) on $[0,T]$ with $1+q$ bounded away from zero and,
$$
\begin{aligned}
&q\in \widetilde{C}_{T}(\dot{B}^\fd_{2,1})\;\;\;\mbox{and}\;\;\;\;u\in \widetilde{C}_{T}(\dot{B}^{\fd-1}_{2,2})
\cap L^{1}_{T}(\dot{B}^{\fd+1}_{2,2})\cap \widetilde{C}_{T}(\dot{B}^{-1}_{\infty,1})
\cap L^{1}_{T}(\dot{B}^{1}_{\infty,1}).
\end{aligned}
$$
In addition ${\rm div}u$ belongs to $\widetilde{L}^{\infty}_{T}(\dot{B}^{\fd-1}_{2,1})\cap L^{1}_{T}(\dot{B}^{\fd+1}_{2,1})$.

Moreover if $P'(1)>0$, there exists a constant $\ee_{0}$ such that if:
$$\|q_{0}\|_{\dot{B}^{\fd-1}_{2,1}\cap \dot{B}^\fd_{2,1}}+\|u_{0}\|_{\dot{B}^{\fd-1}_{2,2}}\leq\ee_{0},$$
then the solution is global.}
\label{theo1}
\end{thm}
\begin{rem}\sl{
We would like to mention that this theorem could easily be extend to the case of Besov space constructed on Lebesgue spaces with index $p\ne 2$. More precisely as in \cite{DL}, we could obtain the existence of solution under the condition than $1\leq p<2N$ and the uniqueness for $1\leq p\leq N$ when $(q_{0},u_{0})$ belong to $\dot{B}^{\frac{d}{p}}_{p,1}\times (\dot{B}^{\frac{d}{p}-1}_{p,2}\cap \dot{B}^{-1}_{\infty,1})$. These restrictions on the size of $p$ are essentially due to some limitations on the use of the paraproduct when we are dealing with the convection term $u\cdot\n u$.}
\end{rem}
\begin{rem}\sl{
In \cite{BD1} the authors obtain a new entropy inequality for the shallow-water system.
They show that we can control $\sqrt{\rho}$ in $L^{\infty}(L^{2})$. Roughly speaking it means that we control the density $\rho$ in $L^{\infty}(\dot{H}^{1})$. It means in particular that our initial data are very close to the energy data in dimension $d=2$. Indeed we essentially assume only additional condition in terms of vacuum ($\frac{1}{\rho_{0}}\in L^{\infty}$) and a control on the $L^{\infty}$ norm of $\rho_{0}$. This conditions are natural in order to deal with the non linear terms but also for preserving the parabolicity of the momentum equation. Moreover we suppose that $u_{0}$ is in $\dot{B}^{-1}_{\infty,1}$. We would like to emphasize that this last condition is quite optimal for Navier-Stokes system. Indeed in the case of the incompressible Navier-Stokes equations Bourgain and Pavlovic have proved that the system is ill-posed when $u_{0}\in \dot{B}^{-1}_{\infty,\infty}$ (see \cite{BP}).}
\end{rem}
\begin{rem}\sl{
This result is also very interesting in the case of the convergence of the solution of compressible system to the incompressible system when the Mach number goes to $0$. Indeed we are able to deal with initial data of Fujita-Kato type, it means here $u_{0}\in \dot{H}^{\N-1}\cap \dot{B}^{-1}_{\infty,1}$ which improves the result of R.Danchin (see \cite{Dmach}) who needs to assume $u_{0}\in \dot{B}^{\N-1}_{2,1}$.}
\end{rem}
\begin{rem}
\sl{The additionnal assumption $\div u_0\in\dot{B}^{\fd-2}_{2,1}$ is quite natural as it expresses the compressibility of the fluid.
} 
\end{rem}

The previous theorem can be easily adapted to the incompressible density dependent Navier-Stokes equations. We recall here the equations:
\begin{equation}
\begin{cases}
\begin{aligned}
&\p_{t}\rho+{\rm div}(\rho u)=0,\\
&\p_{t}(\rho u)+{\rm div}(\rho u\otimes u)-{\rm div}(2\rho Du)+\n\Pi=\rho f,\\
&{\rm div}u=0,\\
&(\rho,u)_{/t=0}=(\rho_{0},u_{0}).
\end{aligned}
\end{cases}
\label{NSC}
\end{equation}
Following the same ideas than in theorem \ref{theo1}, we obtain the following result:
\begin{thm}\sl{
Let $d\geq 3$. Assume that $u_{0}\in \dot{B}^{\fd-1}_{2,2}\cap \dot{B}^{-1}_{\infty,1}$ with $\div u_{0}=0$, $q_{0}\in \dot{B}^{\fd}_{2,1}$ and that there exists $c>0$ such that $\rho_{0}\geq c$. Then there exists a time $T$ such that there exists a unique solution $(q,u)$ for system (\ref{0.1}) on $[0,T]$ with $1+q$ bounded away from zero and,
$$
\begin{aligned}
&q\in \widetilde{C}_{T}(\dot{B}^\fd_{2,1})\;\;\;\mbox{and}\;\;\;\;u\in \widetilde{C}_{T}(\dot{B}^{\fd-1}_{2,2})
\cap L^{1}_{T}(\dot{B}^{\fd+1}_{2,2})\cap \widetilde{C}_{T}(\dot{B}^{-1}_{\infty,1})
\cap L^{1}_{T}(\dot{B}^{1}_{\infty,1})\\
&\hspace{9cm}\mbox{and}\;\; \n\Pi\in
\widetilde{L}^{1}_{T}(\dot{B}^{\frac{d}{2}-1}_{2,2}).
\end{aligned}
$$
Moreover there exists a constant $\ee_{0}$ such that if:
$$\|q_{0}\|_{\dot{B}^{\fd}_{2,1}}+\|u_{0}\|_{\dot{B}^{\fd-1}_{2,2}}\leq\ee_{0},$$
then the solution is global.}
\label{theo2}
\end{thm}
\begin{rem}\sl{
This result may be considered as an extension of the Fujita-Kato theorem (see \cite{FK}) to the case of incompressible density dependent Navier-Stokes equations. In particular it improves the analysis of \cite{AP} and \cite{H} inasmuch as we deal with a velocity in a Besov space such that the third index is different of $1$ and with a critical initial density in terms of scaling.}
\end{rem}
Our paper is structured as follows. In section \ref{section2}, we give a few notation and briefly introduce the basic Fourier analysis
techniques needed to prove our result. In section \ref{section3} and section \ref{section4}, we prove theorem \ref{theo1}  and more particular the existence of such solution in section \ref{section3} and the uniqueness in section \ref{section4}. In section \ref{section5} we are proving the global well-posedness of theorem \ref{theo1}. In section \ref{section6}, we are dealing with theorem \ref{theo2}.

\section{Littlewood-Paley theory and Besov spaces}
\label{section2}
As usual, the Fourier transform of $u$ with respect to the space variable will be denoted by $\mathcal{F}(u)$ or $\hat{u}$. 
In this section we will state classical definitions and properties concerning the homogeneous dyadic decomposition with respect to the Fourier variable. We will recall some classical results and we refer to \cite{BCD} (Chapter 2) for proofs (and more general properties).

To build the Littlewood-Paley decomposition, we need to fix a smooth radial function $\chi$ supported in (for example) the ball $B(0,\frac{4}{3})$, equal to 1 in a neighborhood of $B(0,\frac{3}{4})$ and such that $r\mapsto \chi(r.e_r)$ is nonincreasing over $\R_+$. So that if we define $\varphi(\xi)=\chi(\xi/2)-\chi(\xi)$, then $\varphi$ is compactly supported in the annulus $\{\xi\in \R^d, \frac{3}{4}\leq |\xi|\leq \frac{8}{3}\}$ and we have that,
\begin{equation}
 \forall \xi\in \R^d\setminus\{0\}, \quad \sum_{l\in\Z} \varphi(2^{-l}\xi)=1.
\label{LPxi}
\end{equation}
Then we can define the \textit{dyadic blocks} $(\ddl)_{l\in \Z}$ by $\ddl:= \varphi(2^{-l}D)$ (that is $\hat{\ddl u}=\varphi(2^{-l}\xi)\hat{u}(\xi)$) so that, formally, we have
\begin{equation}
u=\Sum_l \ddl u
\label{LPsomme} 
\end{equation}
As (\ref{LPxi}) is satisfied for $\xi\neq 0$, the previous formal equality holds true for tempered distributions \textit{modulo polynomials}. A way to avoid working modulo polynomials is to consider the set $\cS_h'$ of tempered distributions $u$ such that
$$
\lim_{l\rightarrow -\infty} \|\dot{S}_l u\|_{L^\infty}=0,
$$
where $\dot{S}_l$ stands for the low frequency cut-off defined by $\dot{S}_l:= \chi(2^{-l}D)$. If $u\in \cS_h'$, (\ref{LPsomme}) is true and we can write that $\dot{S}_l u=\Sum_{k\leq l-1} \ddq u$. We can now define the homogeneous Besov spaces used in this article:
\begin{defi}
\label{LPbesov}
 For $s\in\R$ and  
$1\leq p,r\leq\infty,$ we set
$$
\|u\|_{\dot B^s_{p,r}}:=\bigg(\sum_{l} 2^{rls}
\|\Delta_l  u\|^r_{L^p}\bigg)^{\frac{1}{r}}\ \text{ if }\ r<\infty
\quad\text{and}\quad
\|u\|_{\dot B^s_{p,\infty}}:=\sup_{l} 2^{ls}
\|\Delta_l  u\|_{L^p}.
$$
We then define the space $\dot B^s_{p,r}$ as the subset of  distributions $u\in {\cS}'_h$ such that $\|u\|_{\dot B^s_{p,r}}$ is finite.
\end{defi}
Once more, we refer to \cite{BCD} (chapter $2$) for properties of the inhomogeneous and homogeneous Besov spaces. Among these properties, let us mention:
\begin{itemize}
\item for any $p\in[1,\infty]$ we have the following chain of continuous embeddings:
$$
\dot B^0_{p,1}\hookrightarrow L^p\hookrightarrow \dot B^0_{p,\infty};
$$
\item if $p<\infty$ then 
  $B^{\frac dp}_{p,1}$ is an algebra continuously embedded in the set of continuous 
  functions decaying to $0$ at infinity;
    \item for any  smooth homogeneous  of degree $m$ function $F$ on $\R^d\setminus\{0\}$
the operator $F(D)$ maps  $\dot B^s_{p,r}$ in $\dot B^{s-m}_{p,r}.$ This implies that the gradient operator maps $\dot B^s_{p,r}$ in $\dot B^{s-1}_{p,r}.$  
  \end{itemize}
We refer to \cite{BCD} (lemma 2.1) for the Bernstein lemma (describing how derivatives act on spectrally localized functions), that entails the following embedding result:
\begin{prop}\label{LP:embed}
\sl{For all $s\in\R,$ $1\leq p_1\leq p_2\leq\infty$ and $1\leq r_1\leq r_2\leq\infty,$
  the space $\dot B^{s}_{p_1,r_1}$ is continuously embedded in 
  the space $\dot B^{s-d(\frac1{p_1}-\frac1{p_2})}_{p_2,r_2}.$}
\end{prop}
Then we have:
$$
\dot B^\fdp_{p,1}\hookrightarrow \dot B^0_{\infty,1}\hookrightarrow L^\infty.
$$
In this paper, we shall mainly work with functions or distributions depending on both the time variable $t$ and the space variable $x.$ We shall denote by $\cC(I;X)$ the set of continuous functions on $I$ with values in $X.$ For $p\in[1,\infty]$, the notation $L^p(I;X)$ stands for the set of measurable functions on  $I$ with values in $X$ such that $t\mapsto \|f(t)\|_X$ belongs to $L^p(I)$.

In the case where $I=[0,T],$  the space $L^p([0,T];X)$ (resp. $\cC([0,T];X)$) will also be denoted by $L_T^p X$ (resp. $\cC_T X$). Finally, if $I=\R^+$ we shall alternately use the notation $L^p X.$

The Littlewood-Paley decomposition enables us to work with spectrally localized (hence smooth) functions rather than with rough objects. We naturally obtain bounds for each dyadic block in spaces of type $L^\rho_T L^p.$  Going from those type of bounds to estimates in  $L^\rho_T \dot B^s_{p,r}$ requires to perform a summation in $\ell^r(\Z).$ When doing so however, we \emph{do not} bound the $L^\rho_T \dot B^s_{p,r}$ norm for the time integration has been performed \emph{before} the $\ell^r$ summation.
This leads to the following notation (after J.-Y. Chemin and N. Lerner in \cite{CL}):

\begin{defi}\label{d:espacestilde}
For $T>0,$ $s\in\R$ and  $1\leq r,\rho\leq\infty,$
 we set
$$
\|u\|_{\tilde L_T^\rho \dot B^s_{p,r}}:=
\bigl\Vert2^{js}\|\ddq u\|_{L_T^\rho L^p}\bigr\Vert_{\ell^r(\Z)}.
$$
\end{defi}
One can then define the space $\tilde L^\rho_T \dot B^s_{p,r}$ as the set of  tempered distributions $u$ over $(0,T)\times \R^d$ such that $\lim_{q\rightarrow-\infty}\dot S_q u=0$ in $L^\rho([0,T];L^\infty(\R^d))$ and $\|u\|_{\tilde L_T^\rho \dot B^s_{p,r}}<\infty.$ The letter $T$ is omitted for functions defined over $\R^+.$ 
The spaces $\tilde L^\rho_T \dot B^s_{p,r}$ may be compared with the spaces  $L_T^\rho \dot B^s_{p,r}$ through the Minkowski inequality: we have
$$
\|u\|_{\tilde L_T^\rho \dot B^s_{p,r}}
\leq\|u\|_{L_T^\rho \dot B^s_{p,r}}\ \text{ if }\ r\geq\rho\quad\hbox{and}\quad
\|u\|_{\tilde L_T^\rho \dot B^s_{p,r}}\geq
\|u\|_{L_T^\rho \dot B^s_{p,r}}\ \text{ if }\ r\leq\rho.
$$
All the properties of continuity for the product and composition which are true in Besov spaces remain true in the above  spaces. The time exponent just behaves according to H\"older's inequality. 
\medbreak
Let us now recall a few nonlinear estimates in Besov spaces. Formally, any product of two distributions $u$ and $v$ may be decomposed into 
\begin{equation}\label{eq:bony}
uv=T_uv+T_vu+R(u,v), \mbox{ where}
\end{equation}
$$
T_uv:=\sum_l\dot S_{l-1}u\ddl v,\quad
T_vu:=\sum_l \dot S_{l-1}v\ddl u\ \hbox{ and }\ 
R(u,v):=\sum_l\sum_{|l'-l|\leq1}\ddl u\,\dot\Delta_{l'}v.
$$
The above operator $T$ is called ``paraproduct'' whereas $R$ is called ``remainder''. The decomposition \eqref{eq:bony} has been introduced by J.-M. Bony in \cite{BJM}.

In this article we will frequently use the following estimates (we refer to \cite{BCD} section 2.6, \cite{DG}, for general statements, more properties of continuity for the paraproduct and remainder operators, sometimes adapted to $\tilde L_T^\rho \dot B^s_{p,r}$ spaces): under the same assumptions there exists a constant $C>0$ such that if $1/p_1+1/p_2=1/p$, and $1/r_1+1/r_2=1/r$:
$$
\|\dot{T}_u v\|_{\dot{B}_{2,1}^s}\leq C \|u\|_{L^\infty} \|v\|_{\dot{B}_{2,1}^s},
$$
$$\|\dot{T}_u v\|_{\dot{B}_{p,r}^{s+t}}\leq C\|u\|_{\dot{B}_{p_1,r_1}^t} \|v\|_{\dot{B}_{p_2,r_2}^s} \quad (t<0),
$$
\begin{equation}
 \|\dot{R}(u,v)\|_{\dot{B}_{p,r}^{s_1+s_2-\fd}} \leq C\|u\|_{\dot{B}_{p_1,r_1}^{s_1}} \|v\|_{\dot{B}_{p_2,r_2}^{s_2}} \quad (s_1+s_2>0).
\label{estimbesov}
\end{equation}
Let us now turn to the composition estimates. We refer for example to \cite{BCD} (Theorem $2.59$, corollary $2.63$)):
\begin{prop}
\sl{\begin{enumerate}
 \item Let $s>0$, $u\in \dot{B}_{2,1}^s\cap L^{\infty}$ and $F\in W_{loc}^{[s]+2, \infty}(\R^d)$ such that $F(0)=0$. Then $F(u)\in \dot{B}_{2,1}^s$ and there exists a function of one variable $C_0$ only depending on $s$, $d$ and $F$ such that
$$
\|F(u)\|_{\dot{B}_{2,1}^s}\leq C_0(\|u\|_{L^\infty})\|u\|_{\dot{B}_{2,1}^s}.
$$
\item If $u$ and $v\in\dot{B}_{2,1}^\fd$ and if $v-u\in \dot{B}_{2,1}^s$ for $s\in]-\fd, \fd]$ and $G\in W_{loc}^{[s]+3, \infty}(\R^d)$, then $G(v)-G(u)$ belongs to $\dot{B}_{2,1}^s$ and there exists a function of two variables $C$ only depending on $s$, $d$ and $G$ such that
$$
\|G(v)-G(u)\|_{\dot{B}_{2,1}^s}\leq C(\|u\|_{L^\infty}, \|v\|_{L^\infty})\left(|G'(0)| +\|u\|_{\dot{B}_{2,1}^\fd} +\|v\|_{\dot{B}_{2,1}^\fd}\right) \|v-u\|_{\dot{B}_{2,1}^s}.
$$
\end{enumerate}}
\label{estimcompo}
\end{prop}
Let us now recall a result of interpolation which explains the link between the space $B^{s}_{p,1}$ and the space $B^{s}_{p,\infty}$ (see
\cite{DFourier} or \cite{BCD} sections $2.11$ and $10.2.4$):
\begin{prop}\sl{
\label{interpolationlog}
There exists a constant $C$ such that for all $s\in\R$, $\e>0$ and
$1\leq p<+\infty$,
$$\|u\|_{\widetilde{L}_{T}^{\rho}(\dot{B}^{s}_{p,1})}\leq C\frac{1+\e}{\e}\|u\|_{\widetilde{L}_{T}^{\rho}(\dot{B}^{s}_{p,\infty})}
\log\biggl(e+\frac{\|u\|_{\widetilde{L}_{T}^{\rho}(\dot{B}^{s-\e}_{p,\infty})}+ \|u\|_{\widetilde{L}_{T}^{\rho}(\dot{B}^{s+\e}_{p,\infty})}}
{\|u\|_{\widetilde{L}_{T}^{\rho}(\dot{B}^{s}_{p,\infty})}}\biggl).$$ \label{5Yudov}}
\end{prop}
Let us end this section by recalling the following estimates for the heat equation:
\begin{prop}\sl{
\label{chaleur} Let $s\in\R$, $(p,r)\in[1,+\infty]^{2}$ and
$1\leq\rho_{2}\leq\rho_{1}\leq+\infty$. Assume that $u_{0}\in \dot{B}^{s}_{p,r}$ and $f\in\widetilde{L}^{\rho_{2}}_{T}
(\dot{B}^{s-2+2/\rho_{2}}_{p,r})$.
Let u be a solution of:
$$
\begin{cases}
\begin{aligned}
&\p_{t}u-\mu\D u=f\\
&u_{t=0}=u_{0}\,.\\
\end{aligned}
\end{cases}
$$
Then there exists $C>0$ depending only on $N,\mu,\rho_{1}$ and
$\rho_{2}$ such that:
$$\|u\|_{\widetilde{L}^{\rho_{1}}_{T}(\dot{B}^{s+2/\rho_{1}}_{p,r})}\leq C\big(
 \|u_{0}\|_{\dot{B}^{s}_{p,r}}+\mu^{\frac{1}{\rho_{2}}-1}\|f\|_{\widetilde{L}^{\rho_{2}}_{T}
 (\dot{B}^{s-2+2/\rho_{2}}_{p})}\big)\,.$$
 If in addition $r$ is finite then $u$ belongs to $C([0,T],\dot{B}^{s}_{p,r})$.}
\end{prop}

\section{Existence of solution}
\label{section3}
\subsection{A priori estimates}
\label{aprioriqu}
Let us emphasize that we do not use the fact that $P'(1)>0$. Denoting by $G$ the unique primitive of $x\mapsto P'(x)/x$ such that $G(1)=0$, recall that system (\ref{0.6}) now reads:
\begin{equation}
\begin{cases}
\begin{aligned}
&\p_{t}q+u\cdot\n q+(1+q){\rm div}u=0,\\
&\p_{t}u+u\cdot\n u-{\cal A}u-2D(u).\n(\ln(1+q))+\n(G(1+q))=0,
\end{aligned}
\end{cases}
\label{syst}
\end{equation}
where the operator $\mathcal{A}$ is defined by: $\mathcal{A} u =\Delta u +\n \div u$.

Let $u_L$ be the unique global solution of the following linear heat equation:
\begin{equation}
\begin{cases}
\begin{aligned}
&\p_{t}u_L-{\cal A}u_L=0,\\
&{u_L}_{t=0}=u_0,\\
\end{aligned}
\end{cases}
\label{systlin}
\end{equation}
Thanks to the classical heat estimates recalled in (\ref{chaleur}) (we refer for example to \cite{BCD}, lemma $2.4$ and chapter $3$), as $u_0\in \dot{B}_{2,2}^{\fd-1}\cap \dot{B}_{\infty,1}^{\-1}$ we have for all time:
\begin{equation}
u_L\in \Big(\Tilde{L}_t^\infty \dot{B}_{2,2}^{\fd-1} \cap L_t^1 \dot{B}_{2,2}^{\fd+1}\Big) \cap \Big(\Tilde{L}_t^\infty \dot{B}_{\infty,1}^{-1} \cap L_t^1 \dot{B}_{\infty,1}^1\Big), 
\label{reguL1}
\end{equation}
and the corresponding energy estimates. Moreover, as $\div u_0 \in \dot{B}_{2,1}^{\fd-2}$, and as $\div u_L$ satisfies:
\begin{equation}
\begin{cases}
\begin{aligned}
&\p_{t}\div u_L-{\cal A}u_L=\p_{t}\div u_L-2\Delta u_L=0,\\
&{\div u_L}_{|t=0}=\div u_0,\\
\end{aligned}
\end{cases}
\end{equation}
we also have that:
\begin{equation}
\div u_L \in \Tilde{L}_t^\infty \dot{B}_{2,1}^{\fd-2} \cap L_t^1 \dot{B}_{2,1}^{\fd},
\label{reguL2}
\end{equation}
which will be crucial in the study of the density equation. Then, if we denote by $\ub=u-u_L$, we now need to study the following system:
\begin{equation}
\begin{cases}
\begin{aligned}
&\p_t q+(\ub+u_L).\n q +(1+q)\div (\ub+u_L)=0,\\
&\p_t \ub-\Delta \ub-\n \div \ub+ (\ub+u_L).\n \ub+ \ub.\n u_L+u_L.\n u_L & \\
&\hspace{4cm}-2 D(\ub+u_L).\n\Big(\ln(1+q)\Big)+\n \Big(G(1+q)\Big)=0,
\end{aligned}
\end{cases}
\label{systbarre}
\end{equation}
The interest of introducing this system is that the most problematic term in the additionnal external force terms, namely $u_L\cdot \n u_L$, is in fact regular. Thanks to the paraproduct and remainder estimates, we have:
$$
\|u_L\cdot \n u_L\|_{L_t^1 \dot{B}_{2,1}^{\fd-1}} \leq \|u_L\|_{\Tilde{L}_t^\infty \dot{B}_{\infty,1}^{-1}} \|u_L\|_{L_t^1 \dot{B}_{2,2}^{\fd+1}} 
$$ 
In this article we will prove existence and uniqueness of a local solution such that the velocity fluctuation $\ub$ is in the space $\Tilde{L}_T^\infty \dot{B}_{2,1}^{\fd-1} \cap L_T^1 \dot{B}_{2,1}^{\fd+1}$ for some $T>0$. Then, going back to the original functions, thanks to the following embeddings:
$$
\dot{B}_{2,1}^{\fd-1}\hookrightarrow \dot{B}_{2,2}^{\fd-1}, \quad \dot{B}_{2,1}^{\fd-1}\hookrightarrow \dot{B}_{\infty,1}^{-1},
$$
we will end with a velocity $u$ with the same regularity as $u_L$ (that is \ref{reguL1} and \ref{reguL2}).
\\

Let us now state the following transport-diffusion estimates which are adaptations of the ones given in \cite{BCD} section 3:

\begin{lem}\sl{
 Let $T>0$, $-\fd<s<\fd$, $u_0\in\dot{B}_{2,1}^s$, $f\in L_T^1 \dot{B}_{2,1}^s$ and $v,w\in L_T^1 \dot{B}_{2,2}^{\fd+1}\cap \dot{B}_{\infty, 1}^1$. If $u$ is a solution of:
$$
\begin{cases}
\begin{aligned}
&\p_{t}u+v\cdot\n u+u\cdot \n w-{\cal A}u=f,\\
&u_{t=0}=u_0,
\end{aligned}
\end{cases}
$$
then, if $V(t)=\int_0^t \Big(\|\n v(\tau)\|_{\dot{B}_{2,2}^{\fd}\cap \dot{B}_{\infty, 1}^0}+ \|\n w(\tau)\|_{\dot{B}_{2,2}^{\fd}\cap \dot{B}_{\infty, 1}^0}\Big) d\tau$, there exists a constant $C>0$ such that for all $t\in[0,T]$:
$$
\|u\|_{\Tilde{L}_t^\infty \dot{B}_{2,1}^s} +\|u\|_{L_t^1 \dot{B}_{2,1}^{s+2}} \leq C e^{CV(t)} \Big(\|u_0\|_{\dot{B}_{2,1}^s}+ \int_0^t\|f(\tau)\|_{\dot{B}_{2,1}^s} e^{-CV(\tau)} d\tau\Big).
$$}
\end{lem}
\textbf{Proof: } We refer to \cite{BCD} for details. Localizing in frequency, if for $j\in \Z$, $u_j=\ddj u$, we have:
$$
\p_{t}u_j+v\cdot\n u_j-{\cal A}u_j=f_j -\ddj(u\cdot \n w)+R_j,
$$
where $R_j=[v.\n, \ddj]u$. Taking the $L^2$ innerproduct we obtain:
$$
\p_{t}\|u_j\|_{L^2}+2^{2j} \|u_j\|_{L^2}\leq \|\div v\|_{L^\infty} \|u_j\|_{L^2} + \|f_j\|_{L^2} +\|\ddj(u\cdot \n w)\|_{L^2}+\|R_j\|_{L^2}.
$$
Classical commutator estimates (we refer to \cite{BCD} section $2.10$) then imply that there exists a summable positive sequence $c_j=c_j(t)$ whose sum is $1$ such that:
$$
\|R_j\|_{L^2}\leq c_j 2^{-js} \|\n v\|_{\dot{B}_{2,\infty}^\fd \cap \dot{B}_{\infty,1}^0} \|u\|_{\dot{B}_{2,1}^s}.
$$
Thanks to the paraproduct and remainder laws, we have:
$$
\|u\cdot \n w\|_{\dot{B}_{2,1}^s}\leq C\|u\|_{\dot{B}_{2,1}^s} \|\n w\|_{\dot{B}_{2,2}^{\fd}\cap \dot{B}_{\infty, 1}^0}.
$$
so that we finally obtain the result. $\blacksquare$
\\

Concerning the transport equation of the density fluctuation, localization implies:
$$
\p_{t}q_j+(\ub+u_L)\cdot\n q_j=-\ddj\Big((1+q)\div(\ub+u_L)\Big)+R_j,
$$
where $R_j=[(\ub+u_L).\n, \ddj]q$. Using the same method, we get the estimate:
\begin{multline}
 \|q\|_{\Tilde{L}_t^\infty \dot{B}_{2,1}^\fd} \leq \|q_0\|_{\Tilde{L}_t^\infty \dot{B}_{2,1}^\fd} + \int_0^t \Big(\|q\|_{\dot{B}_{2,1}^\fd} \|\n (\ub+u_L)\|_{\dot{B}_{2,2}^{\fd}\cap \dot{B}_{\infty, 1}^0}\\
+(1+\|q\|_{\dot{B}_{2,1}^\fd})(\|\div u_L\|_{\dot{B}_{2,1}^\fd}+ \|\n \ub\|_{\dot{B}_{2,1}^\fd}) \Big) d\tau.
\end{multline}
\begin{rem}
 \sl{Note that here, due to the external force terms for the density, we need to have a $\dot{B}_{2,1}^\fd$-control of $\div u_L$.}
\end{rem}
With a rough majoration,
$$
\|q\|_{\Tilde{L}_t^\infty \dot{B}_{2,1}^\fd} \leq C\Big(\|q_0\|_{\Tilde{L}_t^\infty \dot{B}_{2,1}^\fd} +\int_0^t (1+\|q\|_{\dot{B}_{2,1}^\fd})V'(\tau) d\tau\Big),
$$
where
$$V(t)=\int_0^t \Big(\|\ub\|_{\dot{B}_{2,1}^{\fd+1}}+ \|\n u_L\|_{\dot{B}_{2,2}^\fd \cap \dot{B}_{\infty,1}^0}+\|\div u_L\|_{\dot{B}_{2,1}^\fd}\Big)d\tau,$$
and then thanks to the Gronwall lemma,
$$
\|q\|_{\Tilde{L}_t^\infty \dot{B}_{2,1}^\fd} \leq C e^{CV(t)} \Big(\|q_0\|_{\Tilde{L}_t^\infty \dot{B}_{2,1}^\fd} +\int_0^t V'(\tau) e^{-CV(\tau)} d\tau\Big),
$$
which gives:
\begin{equation}
 \|q\|_{\Tilde{L}_t^\infty \dot{B}_{2,1}^\fd} \leq C e^{CV(t)}\left(1+\|q_0\|_{\dot{B}_{2,1}^\fd}\right)-1,
\label{aprioriq}
\end{equation}
For the velocity, the same localization technique is used, and if we introduce: for $j\in \Z$, $\ub_j=\ddj \ub$, we have:
$$
\p_{t}\ub_j+(\ub+u_L)\cdot\n \ub_j-{\cal A}\ub_j=\ddj f -\ddj(\ub\cdot \n u_L)+R_j,
$$
where
$$
\begin{cases}
 R_j=[(\ub+u_L).\n, \ddj]\ub,\\
f=-u_L.\n u_L+2 D(\ub+u_L).\n\Big(\ln(1+q)\Big)-\n \Big(G(1+q)\Big).
\end{cases}
$$
Taking the $L^2$ innerproduct with $\ub_j$ we obtain:
$$
\p_{t}\|\ub_j\|_{L^2}+2^{2j} \|\ub_j\|_{L^2}\leq \|\div (\ub+u_L)\|_{L^\infty} \|\ub_j\|_{L^2} + \|\ddj f\|_{L^2} +\|\ddj(\ub\cdot \n u_L)\|_{L^2}+\|R_j\|_{L^2}.
$$
After a time integration, the mutliplication by $2^{j(\fd-1)}$ followed by a summation over $j\in \Z$ gives ($\ub(0)=0$ and the commutator is estimated as above):
\begin{multline}
 \|\ub\|_{\Tilde{L}_t^\infty \dot{B}_{2,1}^{\fd-1}}+ \|\ub\|_{L_t^1 \dot{B}_{2,1}^{\fd+1}} \leq\\
\int_0^t \Big(\|f\|_{\dot{B}_{2,1}^{\fd-1}} +\|\ub\cdot \n u_L\|_{\dot{B}_{2,1}^{\fd-1}} +\|\n (\ub+ u_L)\|_{\dot{B}_{2,2}^{\fd}\cap \dot{B}_{\infty,1}^0} \|\ub\|_{\dot{B}_{2,1}^{\fd-1}}\Big) d\tau.
\end{multline}
Classical computations show that:
$$
\|\ub\cdot \n u_L\|_{\dot{B}_{2,1}^{\fd-1}}\leq \|\ub\|_{\dot{B}_{2,1}^{\fd-1}}\|\n u_L\|_{\dot{B}_{2,2}^{\fd}\cap \dot{B}_{\infty,1}^0},
$$
and from the definition of $f$:
$$
\|f\|_{\dot{B}_{2,1}^{\fd-1}} \leq \|u_L\|_{\dot{B}_{\infty,1}^{-1}} \|u_L\|_{\dot{B}_{2,2}^{\fd+1}} +\|\n G(1+q)\|_{\dot{B}_{2,1}^{\fd-1}} +2\|\n\ln(1+q)\cdot D(\ub+u_L)\|_{\dot{B}_{2,1}^{\fd-1}}.
$$
The second term is estimated thanks to the composition lemma (see ???):
$$
\|\n G(1+q)\|_{\dot{B}_{2,1}^{\fd-1}}\leq \|G(1+q)\|_{\dot{B}_{2,1}^{\fd}}\leq C(\|q\|_{L^\infty})\|q\|_{\dot{B}_{2,1}^{\fd}}.
$$
The last term has to be treated carefully. As we cannot rely on the smallness of $q$ (our data can be large), we will follow the same method as R. Danchin in \cite{DW} and thanks to a frequency cut-off we decompose this term into two parts which will both be small. For $m\in\Z$:
$$
\n\ln(1+q)\cdot D(\ub+u_L)=I+II
\mbox{ with: }
\begin{cases}
I=\n\Big(\ln(1+q)-\ln(1+\dot{S}_m q)\Big)\cdot D(\ub+u_L)\\
II=\n\ln(1+\dot{S}_m q)\cdot D(\ub+u_L).
\end{cases}
$$
Thanks to the paraproduct and remainder laws we have,
\begin{multline}
 \|I\|_{\dot{B}_{2,1}^{\fd-1}} \leq \|\ln(1+q)-\ln(1+\dot{S}_m q)\|_{\dot{B}_{2,1}^\fd}\|\n (\ub+u_L)\|_{\dot{B}_{2,2}^{\fd}\cap \dot{B}_{\infty,1}^0}\\
\leq C(\|q\|_{L^\infty})\|q-\dot{S}_m q\|_{\dot{B}_{2,1}^{\fd}} \|\n (\ub+u_L)\|_{\dot{B}_{2,2}^{\fd}\cap \dot{B}_{\infty,1}^0}.
\label{frecutoff}
\end{multline}
We then use the estimate given in \cite{DW}: there exists a constant $C \geq 1$ such that:
$$
\|q-\dot{S}_m q\|_{\dot{B}_{2,1}^{\fd}} \leq \|q_0-\dot{S}_m q_0\|_{\dot{B}_{2,1}^{\fd}} + (1+\|q_0\|_{\dot{B}_{2,1}^{\fd}}) (e^{CV(t)}-1)
$$
If $0<\alpha<1$, the other term is estimated the following way:
\begin{multline}
 \|II\|_{\dot{B}_{2,1}^{\fd-1}} \leq \|\n\ln(1+\dot{S}_m q)\|_{\dot{B}_{2,1}^{\fd+\alpha-1}} \|D(\ub+u_L)\|_{\dot{B}_{2,1}^{\fd-\alpha}}\\
\leq C(\|q\|_{L^\infty}) 2^{m\alpha} \|q\|_{\dot{B}_{2,1}^\fd} \|\ub+u_L\|_{\dot{B}_{2,2}^{\fd+1-\alpha}}.
\end{multline}
And thanks to the Gronwall lemma, we finally obtain that for all $0<\alpha<1$ (for example we can take $\alpha=1/2$) and all $m\in \Z$:
\begin{multline}
 \|\ub\|_{\Tilde{L}_t^\infty \dot{B}_{2,1}^{\fd-1}}+ \|\ub\|_{L_t^1 \dot{B}_{2,1}^{\fd+1}} \leq\\
e^{CV(t)}\Bigg[\|u_L\|_{\Tilde{L}_t^\infty \dot{B}_{\infty,1}^{-1}} \|u_L\|_{L_t^1 \dot{B}_{2,2}^{\fd+1}} +\int_0^t C(\|q\|_{L^\infty}) \Bigg( \|q\|_{\dot{B}_{2,1}^\fd}\Big(1+2^{\alpha m}\|\ub+u_L\|_{\dot{B}_{2,2}^{\fd+1-\alpha}}\Big)\\
+\Big(\|q_0-S_m q_0\|_{\dot{B}_{2,1}^\fd} +(1+\|q_0\|_{\dot{B}_{2,1}^\fd})(e^{CV(\tau)}-1)\Big) \|\ub+ u_L\|_{\dot{B}_{2,2}^{\fd+1}\cap \dot{B}_{\infty,1}^1} \Bigg) d\tau\Bigg],
\label{aprioriu}
\end{multline}
Thanks to the H\"older estimate, we have:
$$
\|\ub+u_L\|_{L_t^1\dot{B}_{2,2}^{\fd+1-\alpha}}\leq t^{\frac{\alpha}{2}}\|\ub+u_L\|_{\Tilde{L}_t^{\frac{1}{1-\frac{\alpha}{2}}} \dot{B}_{2,2}^{\fd+1-\alpha}},
$$
and by interpolation,
$$
 \|f\|_{\Tilde{L}_t^{\frac{1}{1-\frac{\alpha}{2}}} \dot{B}_{2,2}^{\fd+1-\alpha}} \leq \|f\|_{\Tilde{L}_t^\infty\dot{B}_{2,2}^{\fd-1}}^{\frac{\alpha}{2}} \|f\|_{L_t^1\dot{B}_{2,2}^{\fd+1}}^{1-\frac{\alpha}{2}}\leq \frac{\alpha}{2} \|f\|_{\Tilde{L}_t^\infty\dot{B}_{2,2}^{\fd-1}}+ (1-\frac{\alpha}{2}) \|f\|_{L_t^1\dot{B}_{2,2}^{\fd+1}},
$$
so that,
$$
\|\ub+u_L\|_{L_t^1\dot{B}_{2,2}^{\fd+1-\alpha}}\leq t^{\frac{\alpha}{2}} (\|\ub+u_L\|_{\Tilde{L}_t^\infty\dot{B}_{2,2}^{\fd-1}}+\|\ub+u_L\|_{L_t^1\dot{B}_{2,2}^{\fd+1}}),
$$
If we denote $\beta(t)= \|\ub\|_{\Tilde{L}_t^\infty\dot{B}_{2,2}^{\fd-1}}+\|\ub\|_{L_t^1\dot{B}_{2,2}^{\fd+1}}$, the estimates on the velocity turn into:
\begin{multline}
 \beta(t)\leq e^{CV(t)}\Bigg(\|u_L\|_{\Tilde{L}_t^\infty \dot{B}_{\infty,1}^{-1}} \|u_L\|_{L_t^1 \dot{B}_{2,2}^{\fd+1}}\\
+C(\|q\|_{\Tilde{L}_t^\infty \dot{B}_{2,1}^\fd})\Big[t\|q\|_{\Tilde{L}_t^\infty \dot{B}_{2,1}^\fd}+\Big(\|q_0-S_m q_0\|_{\dot{B}_{2,1}^\fd} +(1+\|q_0\|_{\dot{B}_{2,1}^\fd})(e^{CV(t)}-1)+2^{\alpha m} t^{\frac{\alpha}{2}}\|q\|_{\Tilde{L}_t^\infty \dot{B}_{2,1}^\fd} \Big)\\
\times\Big(\beta(t)+ \|u_L\|_{\Tilde{L}_t^\infty\dot{B}_{2,2}^{\fd-1}\cap \dot{B}_{\infty,1}^{-1}} +\|u_L\|_{L_t^1\dot{B}_{2,2}^{\fd+1}\cap \dot{B}_{\infty,1}^1}\Big)\Big]\Bigg).
\end{multline}
Let us now state and prove the following lemma:
\begin{lem}
 \sl{Let $(q,u)$ satisfying $(SW)$ on $[0,T]\times \R^d$. Assume that $q\in\cC^1([0,T], \dot{B}_{2,1}^\fd)$ and $u \mbox{, } u_L\in\cC^1([0,T], \dot{B}_{2,2}^{\fd-1}\cap \dot{B}_{\infty,1}^{-1}\cap \dot{B}_{2,2}^{\fd+1} \cap \dot{B}_{\infty,1}^1))^d$, where $u_L$ satisfies:
\begin{equation}
\d_t u_L-\cA u_L=0, \quad {u_L}_{|t=0}=u_0.
\label{uL}
\end{equation}
If we denote $\ub=u-u_L$, there exist three positive constants $\eta$, $C$ (only depending on $d$), $C'$ (depending on $C$, $d$ and $q_0$) and $m\in \Z$ ($=m(\eta,q_0)$) such that if $\eta\in]0,1]$ satisfies:
\begin{equation}
 \eta e^{3C\eta}(1+C')(1+\|u_0\|_{\dot{B}_{2,2}^{\fd-1}\cap \dot{B}_{\infty,1}^{-1}}) \leq \frac{1}{2},
\end{equation}
and $m$ is chosen such that $\|q_0-S_m q_0\|_{\dot{B}_{2,1}^\fd}\leq \eta^2$. Then if $T$ is small enough so that:
\begin{equation}
 \begin{cases}
  \int_0^T \Big(\|\n u_L\|_{\dot{B}_{2,2}^\fd \cap \dot{B}_{\infty,1}^0}+\|\div u_L\|_{\dot{B}_{2,1}^\fd}\Big)d\tau \leq \eta^2,\\
TC'+(1+C')(e^{CV(T)}-1)+2^{\alpha m} T^{\frac{\alpha}{2}} C' \leq \eta^2.
 \end{cases}
\end{equation}
then we have, for all $t\in [0,T]$,
\begin{equation}
 \begin{cases}
  \|q\|_{\Tilde{L}_T^\infty \dot{B}_{2,1}^\fd}\leq e^{3C\eta}\left(1+\|q_0\|_{\dot{B}_{2,1}^\fd}\right)-1,\\
\|\ub\|_{\Tilde{L}_t^\infty \dot{B}_{2,1}^{\fd-1}}+ \|\ub\|_{L_t^1 \dot{B}_{2,1}^{\fd+1}} \leq 2\eta.
 \end{cases}
\end{equation}
}
\label{estimap}
\end{lem}

\textbf{Proof: } We recall that $(q,\ub)$ satisfies on the interval $[0,T]$ the following system:
$$
\begin{cases}
\begin{aligned}
&\p_t q+(\ub+u_L).\n q +(1+q)\div (\ub+u_L)=0,\\
&\p_t \ub-\Delta \ub-\n \div \ub+ (\ub+u_L).\n \ub+ \ub.\n u_L+u_L.\n u_L & \\
&\hspace{4cm}-2 D(\ub+u_L).\n\Big(\ln(1+q)\Big)+\n \Big(G(1+q)\Big)=0,
\end{aligned}
\end{cases}
$$
Assume that $T$ is small enough so that we have for some $\eta\in [0,1]$ (to be precised later):
$$
\int_0^T \Big(\|\n u_L\|_{\dot{B}_{2,2}^\fd \cap \dot{B}_{\infty,1}^0}+\|\div u_L\|_{\dot{B}_{2,1}^\fd}\Big)d\tau \leq \eta^2\leq \eta,
$$
and let us define
$$
T^*=\Sup \{t\in[0,T], \quad \int_0^t \|\ub\|_{\dot{B}_{2,1}^{\fd+1}}\leq 2\eta\}.
$$
Then for all $t\in[0,T^*[$, from (\ref{aprioriq}) we have ($\eta\leq 1$):
\begin{equation}
 \|q\|_{\Tilde{L}_t^\infty \dot{B}_{2,1}^\fd} \leq e^{C(2\eta + \eta^2)}\left(1+\|q_0\|_{\dot{B}_{2,1}^\fd}\right)-1\leq e^{3C}\left(1+\|q_0\|_{\dot{B}_{2,1}^\fd}\right)-1\overset{\mbox{def}}{=}c_0.
\end{equation}
Concerning the velocity, using the previous estimate on $q$ and the fact that
$$\|u_L\|_{\Tilde{L}_t^\infty\dot{B}_{2,2}^{\fd-1}\cap \dot{B}_{\infty,1}^{-1}} +\|u_L\|_{L_t^1\dot{B}_{2,2}^{\fd+1}\cap \dot{B}_{\infty,1}^{1}} \leq \|u_0\|_{\dot{B}_{2,2}^{\fd-1}\cap\dot{B}_{\infty,1}^{-1}}$$
allows us to write that:
\begin{multline}
 \beta(t)\leq e^{3C\eta}\Bigg(\|u_0\|_{\dot{B}_{\infty,1}^{-1}} \eta^2 +C(c_0)\Big[tc_0+\Big(\|q_0-S_m q_0\|_{\dot{B}_{2,1}^\fd} +(1+c_0)(e^{CV(t)}-1)+2^{\alpha m} t^{\frac{\alpha}{2}} c_0 \Big)\\
\times\Big(\beta(t)+ \|u_0\|_{\dot{B}_{2,2}^{\fd-1}\cap \dot{B}_{\infty,1}^{-1}}\Big)\Big]\Bigg).
\end{multline}
Let us first fix $m=m(\eta)$ such that $\|q_0-S_m q_0\|_{\dot{B}_{2,1}^\fd}\leq \eta^2$. Then let us take $T$ so small that:
$$
Tc_0+(1+c_0)(e^{CV(T)}-1)+2^{\alpha m} T^{\frac{\alpha}{2}} c_0 \leq \eta^2.
$$
The estimate turns into:
$$
\beta(t)\leq e^{3C\eta}\left(\|u_0\|_{\dot{B}_{\infty,1}^{-1}} \eta^2 +C(c_0)\eta^2\Big[1+\beta(t) +\|u_0\|_{\dot{B}_{2,2}^{\fd-1}\cap \dot{B}_{\infty,1}^{-1}}\Big] \right).
$$
Then, if $\eta$ is taken so small that
$$\eta e^{3C\eta}(\|u_0\|_{\dot{B}_{\infty,1}^{-1}}+C(c_0)(1+\|u_0\|_{\dot{B}_{2,2}^{\fd-1}\cap \dot{B}_{\infty,1}^{-1}}) \leq \frac{1}{2},$$
then as $\eta\in]0,1]$, we obtain that:
$$
\beta(t)\leq \frac{1}{2}(\beta(t)+\eta),
$$
which implies:
$$
\beta(t)\leq \eta.
$$
Then by contradiction this leads to $T^*=T$ and for all $t\in[0,T]$,
$$
 \|u\|_{\Tilde{L}_t^\infty \dot{B}_{2,1}^{\fd-1}}+ \|u\|_{L_t^1 \dot{B}_{2,1}^{\fd+1}}\leq 2\eta.
$$

\subsection{Existence}
We use a standard scheme for proving the existence of the solutions:
\begin{enumerate}
\item We smooth out the data and get a sequence of smooth solutions $(q^{n},u^{n})_{n\in\mathbb{N}}$ to (\ref{0.6})
on a bounded interval $[0,T^{n}]$ which may depend on $n$. 
\item We exhibit a positive lower bound $T$ for $T^{n}$, and prove uniform estimates on $(q^{n},\bar{u}^{n})$ (we refer to the next subsection for the definition of $\bar{u}^{n}$) in the space
\begin{equation}
 E_{T}=\widetilde{C}_{T}(\dot{B}^{\fd}_{2,1})\times\big(\widetilde{C}_{T}(\dot{B}^\fd_{2,1})\cap L^1_T(\dot{B}^{\fd+1}_{2,1})\big).
\label{ET}
\end{equation}
\item We use compactness to prove that the sequence $(q^{n},u^{n})$ converges, up to extraction, to a solution of (\ref{0.6}).
\end{enumerate}
\subsubsection{Step 1: Friedrichs approximation}

In order to construct approximated solutions of system (\ref{0.6}) we shall use the classical Friedrichs approximation where we define the frequency truncation operator $J_n$ by: 
$$\mbox{for all}\; n\in \mathbb{N}\;\;\mbox{and for all}\;g\in L^2(\R^d),
J_n g= \mathcal{F}^{-1}\left(\textbf{1}_{\frac{1}{n}\leq |\xi|\leq n}(\xi) \hat{g}(\xi)\right),
$$
and we define the following approximated system:
$$
\begin{cases}
\begin{aligned}
&\p_t q_n+J_n(J_n u_n\cdot\n J_n q_n)+J_n\Big((1+J_n q_n){\rm div} J_n u_n\Big)=0,\\
&\p_t u_n+J_n(J_n u_n\cdot\n J_n u_n)-{\cal A}J_n u_n-2J_n\Big(D(J_n u_n).\n(\ln(1+J_n q_n))\Big)\\
&\hspace{8cm} +J_n\Big(\n(G(1+J_n q_n))\Big)=0,\\
&(q_{n},u_{n})_{t=0}=(J_{n}q_{0},J_{n}u_{0}).
\end{aligned}
\end{cases}
$$
we recall the operator $\mathcal{A}$ is defined by: $\mathcal{A} u =\Delta u +\n \div u$.

We can easily check that it is an ordinary differential equation in $L_n^2\times(L_n^2)^d$, where $L_n^2=\{u\in L^2(\R^d), \mbox{ } J_n u=u\}$. Then for every $n\in \mathbb{N}$, by Cauchy-Lipschitz theorem there exists a unique maximal solution in the space $\cC^1([0, T_n^*[, L_n^2)$ and this system can be rewritten into:
\begin{equation}
\begin{cases}
\begin{aligned}
&\p_t q_n+J_n(u_n\cdot\n q_n)+J_n\Big((1+ q_n){\rm div} u_n\Big)=0,\\
&\p_t u_n+J_n(u_n\cdot\n u_n)-{\cal A}u_n-2J_n\Big(D(u_n).\n(\ln(1+q_n))\Big)\\
&\hspace{8cm} +J_n\Big(\n(G(1+q_n))\Big)=0.
\end{aligned}
\end{cases}
\label{linea}
\end{equation}

\subsubsection{Step 2: Uniform estimates}
In the sequel, we will split $u_{n}$ into the solution of a linear system with initial data $J_{n}u_{0}$,
and the discrepancy to that solution. More precisely, we define by $u^{n}_{L}$ the solution
of the following heat equation:
\begin{equation}
\begin{aligned}
&\p_{t}u^{n}_{L}-{\cal{A}}u^{n}_{L}=0\\
&(u^{n}_{L})_{/t=0}=J_{n}u_{0}.
\end{aligned}
\label{chaleurlineaire}
\end{equation}
We now set $\bar{u}_{n}=u_{n}-u^{n}_{L}$. Obviously, the definition of $\bar{u}_{n}$ leads to the following system:

\begin{equation}
\begin{cases}
\begin{aligned}
&\p_t q_n+J_n\Big((\ub_n+u_L^n).\n q_n\Big) +J_n\Big((1+q_n)\div (\ub_n+u_L^n)\Big)=0,\\
&\p_t \ub_n-\Delta \ub_n-\n \div \ub_n+ J_n\Big((\ub_n+u_L^n).\n \ub_n\Big)+ J_n\Big(\ub_n.\n u_L^n\Big)+J_n\Big(u_L^n.\n u_L^n\Big) & \\
&\hspace{4cm}-2 J_n\Big(D(\ub_n+u_L^n).\n\ln(1+q_n)\Big)+\n J_n\Big(G(1+q_n)\Big)=0,\\
&(q_{n},\bun)_{t=0}=(J_{n}q_{0},0).
\end{aligned}
\end{cases}
\label{linea}
\end{equation}
We would like to obtain uniform estimates on $(q_{n},\bun)$ in the space $E_{T}$ (see \ref{ET}).
Before doing this, let us recall that thanks to proposition \ref{chaleur} and as $J_{n}u_{0}$ uniformly belongs (for all $n$) to $\dot{B}^{-1}_{\infty,1}\cap \dot{B}^{\N-1}_{2,2}$ we obtain that for all $T>0$:
\begin{multline}
\|u^{n}_{L}\|_{\widetilde{L}_{T}^{\infty}(\dot{B}^{\N-1}_{2,2}\cap \dot{B}^{-1}_{\infty,1})}+\|u^{n}_{L}\|_{L^{1}_{T}(\dot{B}^{\N+1}_{2,2}\cap \dot{B}^{1}_{\infty,1})}\\
\leq C(\|u_{0}\|_{\dot{B}^{\N-1}_{2,2}\cap \dot{B}^{-1}_{\infty,1}}+\|f\|_{L_{T}^{1}(\dot{B}^{\N-1}_{2,2}\cap \dot{B}^{-1}_{\infty,1})}).
\label{lineaire1}
\end{multline}
In particular by Besov embedding, we can remark that $\n u^{n}_{L}$ belongs in $L^{1}_{T}(L^{\infty})$, this property will be crucial in the sequel in order to estimate $\bun$. We would also point out that as $\div J_{n}u_{0}$ uniformly belongs (for all $n$) to $\dot{B}^{\fd-1}_{2,1}$ we obtain that for all $T>0$:
\begin{equation}
\begin{aligned}
&\|\div \,u^{n}_{L}\|_{\widetilde{L}_{T}^{\infty}(\dot{B}^{\N-2}_{2,1})} +\|\div\,u^{n}_{L}\|_{L^{1}_{T}(\dot{B}^{\N+1}_{2,1})}\leq C(\|u_{0}\|_{\dot{B}^{\N-1}_{2,2}\cap \dot{B}^{-1}_{\infty,1}}+\|f\|_{L_{T}^{1}(\dot{B}^{\N-1}_{2,2}\cap \dot{B}^{-1}_{\infty,1})}).
\end{aligned}
\label{lineaire2}
\end{equation}
Thanks to the apriori estimates from lemma \ref{estimap} (as $J_n q_n=q_n$ and $J_n \ub_n=\ub_n$ the proof of this lemma, which is based on $L^2$-scalar products, remains true) there exists $\eta>0$ and a time $T>0$ (all of them independant of $n$) such that for any $n\in \mathbb{N}$ and any $t\in[0, \min(T_n^*, T)]$,
\begin{equation}
 \begin{cases}
  \|q_n\|_{\Tilde{L}_T^\infty \dot{B}_{2,1}^\fd}\leq e^{3C\eta}\left(1+\|q_0\|_{\dot{B}_{2,1}^\fd}\right)-1,\\
\|\ub_n\|_{\Tilde{L}_t^\infty \dot{B}_{2,1}^{\fd-1}}+ \|\ub_n\|_{L_t^1 \dot{B}_{2,1}^{\fd+1}} \leq 2\eta.
 \end{cases}
\end{equation}
From this we deduce that the $L_n^2$-norm of $(q_n, \ub_n)$ is bounded. As $J_n$ is the truncation operator in $\{\xi\in\R^d, \frac{1}{n} \leq |\xi|\leq n\}$ the bound blows up as $n$ goes to infinity, but all that is important is that it implies (by contradiction) that for all $n$, the maximal lifespan $T_n^*\geq T$.

\subsubsection{Step 3: Time derivatives}

Once the uniform time $T$ is obtained the rest of the method is very classical. Using the previous uniform estimates to bound the time derivatives of the approximated solutions, we obtain that:
\begin{lem}
 \sl{With the same notations, $(\p_t q_n)_n$ is (uniformly in $n$) bounded in $L_T^2 \dot{B}_{2,1}^{\fd-1}$ and $(\p_t \ub_n)_n$ is bounded in $L_T^\frac{4}{3} \dot{B}_{2,1}^{\fd-\frac{3}{2}}+L_T^\infty \dot{B}_{2,1}^{\fd-1}$ and then in $L_T^\frac{4}{3} (\dot{B}_{2,1}^{\fd-1} +\dot{B}_{2,1}^{\fd-\frac{3}{2}})$.
}
\end{lem}
This result allows to get that:
\begin{itemize}
 \item $q_n-q_n(0)$ is (uniformly in $n$) bounded in $\mathcal{C}_T \dot{B}_{2,1}^\fd\cap \mathcal{C}_T^{\frac{1}{2}} \dot{B}_{2,1}^{\fd-1}$,
\item $\ub_n$ is bounded in $\mathcal{C}_T \dot{B}_{2,1}^{\fd-1}\cap \mathcal{C}_T^{\frac{1}{4}} (\dot{B}_{2,1}^{\fd-1}+\dot{B}_{2,1}^{\fd-\frac{3}{2}})$.
\end{itemize}

\subsubsection{Step 4: compactness and convergence}

This part is also classical and we refer for example to \cite{DW} (chapter $10$) for details: using the previous result and the Ascoli theorem, we can extract a subsequence that weakly converges towards some couple $(q,\ub)$, which is proved to be a solution of the original system and to satisfy the energy estimates. This concludes the existence part of the theorem.

\section{Uniqueness}
\label{section4}

Once more, system $(SW)$ is very close to $(NSC)$ and the uniqueness is dealt the same way except, obviously, that here the external force terms that have to be plugged into the apriori estimates are different from the ones in $(NSC)$ and we will focus on it in this section. As for $(NCS)$ we will have (due to endpoints in the paradifferential remainder) to treat separately the cases $d=2$ and $d\geq 3$. The second difficulty is that, as we present a local result for large data, we will have to make frequency cut-off (as in \cite{BCD}) in order to bound some external force terms.\\

Let us introduce for $s\in \R$ the following space:
$$
E_d(t)=\Tilde{L}_t^\infty(\dot{B}_{2,2}^{\fd-1}\cap \dot{B}_{\infty,1}^{-1}) \cap L_t^1(\dot{B}_{2,2}^{\fd+1}\cap \dot{B}_{\infty,1}^{1}).
$$

\begin{thm}
 \sl{Let $d\geq 2$ and assume that $(q_i, u_i)$ ($i\in\{1,2\}$) are two solutions of $(SW)$ with the same initial data on the same interval $[0, T]$ and both belonging to the space $E_d(T)$. Then $(q_1,u_1)\equiv(q_2,u_2)$ on $[0,T]$.
}
\end{thm}

\textbf{Proof: } for $i\in\{1,2\}$, let us introduce $\ub_i=u_i-u_L$ (see (\ref{uL}) for the definition of $u_L$.), then $(q_i, \ub_i)$ satisfy the system:
$$
\begin{cases}
\begin{aligned}
&\p_t q_i+(\ub_i+u_L).\n q_i +(1+q_i)\div (\ub_i+u_L)=0,\\
&\p_t \ub_i-\Delta \ub_i-\n \div \ub_i+ (\ub_i+u_L).\n \ub_i+ \ub_i.\n u_L+u_L.\n u_L & \\
&\hspace{4cm}-2 D(\ub_i+u_L).\n\Big(\ln(1+q_i)\Big)+\n \Big(G(1+q_i)\Big)=0,
\end{aligned}
\end{cases}
$$
and if we denote by $\dq=q_1-q_2$ and $\dub=\ub_1-\ub_2$, then $(\dq,\dub)$ satisfy the following system:
\begin{equation}
\begin{cases}
\begin{aligned}
&\p_t \dq+(u_L+\ub_2).\n \dq= \de F_1 +\de F_2 +\de F_3,\\
&\p_t \dub-\Delta \dub-\n \div \dub+ (u_L+\ub_1).\n \dub+ \dub.\n (u_L+\ub_2)=\de G_1 +\de G_2 +\de G_3,
\end{aligned}
\end{cases}
\label{SWd}
\end{equation}
with:
$$
\begin{cases}
\de F_1=-\dub.\nabla q_1,\\
\de F_2=-(1+q_1)\div \dub,\\
\de F_3=-\dq.\div(\ub_2+u_L),\\
\de G_1=2 D(u_L+\ub_1).\n\Big(\ln(1+q_1)-\ln(1+q_2)\Big),\\
\de G_2=-D(\dub).\n \ln(1+q_2),\\
\de G_3=-\n \Big(G(1+q_1)-G(1+q_2)\Big).
\end{cases}
$$
\subsection{The case $d\geq 3$}
We wish to prove (as for $(NSC)$) the uniqueness in the following space:
$$
F_T=\cC_T\dot{B}_{2,1}^{\fd-1}\times \Big(\cC_T\dot{B}_{2,1}^{\fd-2}\cap L_T^1\dot{B}_{2,1}^\fd \Big)^d
$$
Due to endpoint estimates for the paradifferential remainder, the case $d=2$ has to be treated in a different space. We refer to the following section. As for the classical $(NSC)$ system, we prove that $(\dq,\dub)\in F_T$ (the proof is left to the reader and the computations are the same as the ones done in the following).
\\

Let us begin with $\dq$. As $q_1$ and $q_2$ have the same initial data, doing the same computations as for the transport estimates (see \cite{BCD} theorem $3.14$) leads to:
\begin{multline}
 \|\dq\|_{\Tilde{L}_t^\infty \dot{B}_{2,1}^{\fd-1}}\\
\leq \int_0^t \left(\Big(\|\div (u_L+\ub_2)\|_{L^\infty}+\|\n (\ub_2+u_L)\|_{\dot{B}_{2,\infty}^\fd \cap L^{\infty}}\Big)\|\dq\|_{\Tilde{L}_\tau^\infty \dot{B}_{2,1}^{\fd-1}} +\|\de F\|_{\dot{B}_{2,1}^{\fd-1}}\right)d\tau,
\end{multline}
where $\de F=\de F_1 +\de F_2 +\de F_3$. These terms are estimated thanks to the paraproduct and remainder estimates recalled in section \ref{section2} (see (\ref{estimbesov})):
\begin{itemize}
 \item Thanks to the Bernstein lemma we have $\dot{B}_{1,1}^{d-1} \hookrightarrow \dot{B}_{2,1}^{\fd-1}$ so that:
\begin{multline}
 \|\de F_1\|_{\dot{B}_{2,1}^{\fd-1}} \leq \|T_{\dub}\n q_1\|_{\dot{B}_{2,1}^{\fd-1}} +\|T_{\n q_1}\dub\|_{\dot{B}_{2,1}^{\fd-1}} +\|R(\dub,\n q_1)\|_{\dot{B}_{1,1}^{d-1}},\\
\leq\|\dub\|_{L^\infty} \|\n q_1\|_{\dot{B}_{2,1}^{\fd-1}} + \|\n q_1\|_{\dot{B}_{\infty,\infty}^{-1}} \|\dub\|_{\dot{B}_{2,1}^\fd}+ \|\dub\|_{\dot{B}_{2,2}^\fd} \|\n q_1\|_{\dot{B}_{2,2}^{\fd-1}}\\
\leq \left(\|\dub\|_{\dot{B}_{\infty,1}^0} +\|\dub\|_{\dot{B}_{2,1}^\fd}\right)\|q_1\|_{\dot{B}_{2,1}^\fd} \leq \|\dub\|_{\dot{B}_{2,1}^\fd} \|q_1\|_{\dot{B}_{2,1}^\fd}.
\end{multline}
\item Similarly we get that
$$
\begin{cases}
 \|\de F_2\|_{\dot{B}_{2,1}^{\fd-1}} \leq \Big(1+\|q_1\|_{\dot{B}_{2,1}^\fd}\Big)\|\dub\|_{\dot{B}_{2,1}^\fd},\\
 \|\de F_3\|_{\dot{B}_{2,1}^{\fd-1}} \leq \|\dq\|_{\dot{B}_{2,1}^{\fd-1}} \|u_L+\ub_2\|_{\dot{B}_{2,2}^{\fd+1} \cap \dot{B}_{\infty,1}^1}
\end{cases},
$$
\end{itemize}
then we obtain:
\begin{multline}
 \|\dq\|_{\Tilde{L}_t^\infty \dot{B}_{2,1}^{\fd-1}}\\
\leq \int_0^t \left(\|u_L+\ub_2\|_{\dot{B}_{2,\infty}^{\fd+1} \cap \dot{B}_{\infty,1}^1}\|\dq\|_{\Tilde{L}_\tau^\infty \dot{B}_{2,1}^{\fd-1}} +\Big(1+\|q_1\|_{\dot{B}_{2,1}^\fd}\Big)\|\dub\|_{\dot{B}_{2,1}^\fd}\right)d\tau.
\end{multline}
Finally, thanks to the Gronwall estimate:
\begin{equation}
 \|\dq\|_{\Tilde{L}_t^\infty \dot{B}_{2,1}^{\fd-1}} \leq e^{C\int_0^t \|\ub_2+u_L\|_{\dot{B}_{2,\infty}^{\fd+1} \cap \dot{B}_{\infty,1}^1}d\tau} \int_0^t \Big(1+\|q_1\|_{\dot{B}_{2,1}^\fd}\Big)\|\dub\|_{\dot{B}_{2,1}^\fd} d\tau.
\label{estimq}
\end{equation}
\begin{rem}
 \sl{Note that this estimate on $\dq$ is valid for all $d\geq 2$.}
\end{rem}

Concerning the velocity, using the a priori estimate for the transport-diffusion equation provided in the present article, we can write that (we recall that $\ub_1$ and $\ub_2$ have the same initial data):
\begin{multline}
\|\dub\|_{\Tilde{L}_t^\infty \dot{B}_{2,1}^{\fd-2}}+ \|\dub\|_{L_t^1 \dot{B}_{2,1}^\fd} \leq e^{C\int_0^t \Big(\|\n \ub_1+\n u_L\|_{\dot{B}_{2,2}^\fd \cap \dot{B}_{\infty,1}^0} +\|\n \ub_2+\n u_L\|_{\dot{B}_{2,2}^\fd \cap \dot{B}_{\infty,1}^0}\Big) d\tau}\\
\times\int_0^t \|\de G_1+\de G_2+ \de G_3\|_{\dot{B}_{2,1}^{\fd-2}} d\tau.
\label{estimu1}
\end{multline}
\begin{itemize}
 \item The last term is dealt the usual way:
\begin{multline}
  \|\de G_3\|_{\dot{B}_{2,1}^{\fd-2}}\leq C(\|q_1\|_{L^\infty}, \|q_2\|_{L^\infty})\Big(1+\|q_1\|_{\dot{B}_{2,1}^\fd}+\|q_2\|_{\dot{B}_{2,1}^\fd}\Big)\|\dq\|_{\dot{B}_{2,1}^{\fd-1}}\\
\leq C(\|q_0\|_{\dot{B}_{2,1}^\fd})\|\dq\|_{\dot{B}_{2,1}^{\fd-1}}.
\end{multline}

\item Without surprise, the first term is estimated by:
\begin{multline}
\|\de G_1\|_{\dot{B}_{2,1}^{\fd-2}}\leq 2 \|D(u_L+\ub_1)\|_{\dot{B}_{2,2}^\fd \cap \dot{B}_{\infty,1}^0}\|\n\Big(\ln(1+q_1)-\ln(1+q_2)\Big)\|_{\dot{B}_{2,1}^{\fd-2}},\\
\leq C(\|q_0\|_{\dot{B}_{2,1}^\fd}) \|u_L+\ub_1\|_{\dot{B}_{2,2}^{\fd+1} \cap \dot{B}_{\infty,1}^1} \|\dq\|_{\dot{B}_{2,1}^{\fd-1}}.
\end{multline}
Note that here, after using the injection $\dot{B}_{1,1}^{d-2} \hookrightarrow \dot{B}_{2,1}^{\fd-2}$, we needed $d-2>0$ in the remainder. In the case $d=2$, this term will have to be dealt differently in the following subsection.
\item As in (\ref{frecutoff}), we have to decompose $\de G_2$ into two parts:
$$
\de G_2=-D(\dub).\n \Big(\ln(1+q_2)-\ln(1+S_m q_2)\Big)-D(\dub).\n \ln(1+S_m q_2)=R_1+R_2.
$$
Here the second paraproduct for $R_2$ requires that $-1+\alpha<0$. The remainders require that $d-2>0$). We obtain that:
\begin{multline}
 \|R_1\|_{\dot{B}_{2,1}^{\fd-2}}\leq \|\dub\|_{\dot{B}_{2,1}^\fd} C(\|q_2\|_{L^\infty})\Big(1+\|q_2\|_{\dot{B}_{2,1}^\fd}\Big)\|q_2-S_m q_2\|_{\dot{B}_{2,1}^\fd}\\
\leq C(\|q_0\|_{\dot{B}_{2,1}^\fd})\|q_2-S_m q_2\|_{\dot{B}_{2,1}^\fd}\|\dub\|_{\dot{B}_{2,1}^\fd}\\
\leq C(\|q_0\|_{\dot{B}_{2,1}^\fd})\Big(\|q_0-S_m q_0\|_{\dot{B}_{2,1}^\fd} +(e^{CV(t)}-1)\Big)\|\dub\|_{\dot{B}_{2,1}^\fd},
\end{multline}
and
\begin{multline}
 \|R_2\|_{\dot{B}_{2,1}^{\fd-2}}\leq \|\dub\|_{\dot{B}_{2,1}^{\fd-\alpha}} \|\ln(1+S_m q_2)\|_{\dot{B}_{2,1}^{\fd+\alpha}}\leq \|\dub\|_{\dot{B}_{2,1}^{\fd-\alpha}} C(\|S_m q_2\|_{L^\infty}) \|S_m q_2\|_{\dot{B}_{2,1}^{\fd+\alpha}}\\
\leq C(\|q_0\|_{\dot{B}_{2,1}^\fd}) 2^{\alpha m}\|\dub\|_{\dot{B}_{2,1}^{\fd-\alpha}}.
\end{multline}
The first term is small if $m$ is large enough and $T$ small enough, and the second term introduces a nonnegative power of $t$.
\end{itemize}
Finally we have:
\begin{multline}
 \|\dub\|_{\Tilde{L}_t^\infty \dot{B}_{2,1}^{\fd-2}}+ \|\dub\|_{L_t^1 \dot{B}_{2,1}^\fd} \leq e^{C\int_0^t \Big(\|\n \ub_1\|_{\dot{B}_{2,1}^\fd} +\|\n \ub_2\|_{\dot{B}_{2,1}^\fd} +\|\n u_L\|_{\dot{B}_{2,2}^\fd \cap \dot{B}_{\infty,1}^0}\Big) d\tau}\\
\times C(\|q_0\|_{\dot{B}_{2,1}^\fd}) \int_0^t \Bigg((1+\|\n \ub_1+\n u_L\|_{\dot{B}_{2,2}^\fd \cap \dot{B}_{\infty,1}^0})\|\dq\|_{\dot{B}_{2,1}^{\fd-1}} +2^{\alpha m}\|\dub\|_{\dot{B}_{2,1}^{\fd-\alpha}}\\
+\Big(\|q_0-S_m q_0\|_{\dot{B}_{2,1}^\fd} +(e^{C\int_0^t \Big(\|\n \ub_2\|_{\dot{B}_{2,1}^\fd} +\|\n u_L\|_{\dot{B}_{2,2}^\fd \cap \dot{B}_{\infty,1}^0}\Big) d\tau}-1)\Big)\|\dub\|_{\dot{B}_{2,1}^\fd}\Bigg)d\tau.
\end{multline}
Introducing $c_0$ a constant only depending on $\|q_0\|_{\dot{B}_{2,1}^\fd}$ and:
$$
V(t)=\int_0^t \Big(\|\n \ub_1\|_{\dot{B}_{2,1}^\fd} +\|\n \ub_2\|_{\dot{B}_{2,1}^\fd} +\|\n u_L\|_{\dot{B}_{2,2}^\fd \cap \dot{B}_{\infty,1}^0}\Big) d\tau.
$$
\begin{multline}
 \|\dub\|_{\Tilde{L}_t^\infty \dot{B}_{2,1}^{\fd-2}}+ \|\dub\|_{L_t^1 \dot{B}_{2,1}^\fd} \leq c_0 e^{CV(t)} \Bigg((t+V(t))\|\dq\|_{\Tilde{L}_t^\infty \dot{B}_{2,1}^{\fd-1}}\\
+2^{\alpha m}\|\dub\|_{L_t^1 \dot{B}_{2,1}^{\fd-\alpha}} +\Big(\|q_0-S_m q_0\|_{\dot{B}_{2,1}^\fd} +(e^{CV(t)}-1)\Big)\|\dub\|_{L_t^1 \dot{B}_{2,1}^\fd}\Bigg).
\end{multline}
If $\beta(t)= \|\ub\|_{\Tilde{L}_t^\infty\dot{B}_{2,2}^{\fd-1}}+\|\ub\|_{L_t^1\dot{B}_{2,2}^{\fd+1}}$, as in section \ref{aprioriqu} we have $\|\dub\|_{L_t^1\dot{B}_{2,2}^{\fd-\alpha}}\leq t^{\frac{\alpha}{2}} \beta(t)$, then as,
\begin{equation}
 \|\dq\|_{\Tilde{L}_t^\infty \dot{B}_{2,1}^{\fd-1}} \leq c_0 e^{CV(t)} \|\dub\|_{L_t^1\dot{B}_{2,1}^\fd},
\label{estimq}
\end{equation}
 we obtain:
$$
\beta(t) \leq c_0 e^{2CV(t)} \Big(t+V(t)+2^{\alpha m}t^{\frac{\alpha}{2}} +\|q_0-S_m q_0\|_{\dot{B}_{2,1}^\fd} +(e^{CV(t)}-1)\Big)\beta(t).
$$
When $\eta\in]0,1]$ satisfies:
\begin{equation}
 2c_0 \eta e^{2C\eta} \leq \frac{1}{2},
\end{equation}
and $m$ is chosen such that $\|q_0-S_m q_0\|_{\dot{B}_{2,1}^\fd}\leq \eta$. Then if $T$ is small enough so that:
\begin{equation}
 \begin{cases}
  \int_0^T \Big(\|\n u_L\|_{\dot{B}_{2,2}^\fd \cap \dot{B}_{\infty,1}^0}+\|\div u_L\|_{\dot{B}_{2,1}^\fd}\Big)d\tau \leq \eta,\\
T+V(T)+(e^{CV(T)}-1)+2^{\alpha m} T^{\frac{\alpha}{2}} \leq \eta.
 \end{cases}
\end{equation}
then we have, for all $t\in [0,T]$,
$$
0\leq \beta(t)\leq 2c_0 \eta e^{2C\eta}\beta(t)\leq \frac{\beta(t)}{2}.
$$
So $\beta(t)=0$ for all $t\in [0,T]$ and the same goes for $\dq$, which proves the uniqueness on $[0,T]$.

\begin{rem}
 \sl{Note that these conditions are implied by the ones from the apriori estimates.}
\end{rem}
To end the proof when $T$ is not small, let us introduce (as in \cite{DW}, section $10.2.4$) the set:
$$
I\overset{\mbox{def}}{=}\{t\in[0,T]/(q_1(t'), \ub_1(t'))=(q_2(t'), \ub_2(t')),\;\forall t'\in[0,t]\}.
$$
This is a nonempty closed subset of $[0,T]$. Using the same method as above allows to prove it is also open and then $I=[0,T]$.

\subsection{The case $d=2$}

In this case, the estimates on $\dq$ remain correct, but the paradifferential remainders, when estimating the external forces terms in the velocity equation, are modified. Indeed, in the case $d=2$ we reach the following endpoint where for all $1/p_1+1/p_2=1=1/r_1+1/r_2$:
$$
\|R(f,g)\|_{\dot{B}_{1,\infty}^0}\leq C\|f\|_{\dot{B}_{p_1,r_1}^{s}}\|g\|_{\dot{B}_{p_2,r_2}^{-s}}.
$$
Estimate (\ref{estimu1}) is then replaced by
\begin{multline}
\|\dub\|_{\Tilde{L}_t^\infty \dot{B}_{2,\infty}^{-1}}+ \|\dub\|_{L_t^1 \dot{B}_{2,\infty}^1} \leq e^{C\int_0^t \Big(\|\n \ub_1+\n u_L\|_{\dot{B}_{2,2}^1 \cap \dot{B}_{\infty,1}^0} +\|\n \ub_2+\n u_L\|_{\dot{B}_{2,2}^1 \cap \dot{B}_{\infty,1}^0}\Big) d\tau}\\
\times\int_0^t \|\de G_1+\de G_2+ \de G_3\|_{\dot{B}_{2,\infty}^{-1}} d\tau,
\end{multline}
with
\begin{multline}
  \|\de G_3\|_{\dot{B}_{2,\infty}^{-1}}\leq C(\|q_1\|_{L^\infty}, \|q_2\|_{L^\infty})\Big(1+\|q_1\|_{\dot{B}_{2,1}^1}+\|q_2\|_{\dot{B}_{2,1}^1}\Big)\|\dq\|_{\dot{B}_{2,\infty}^0}\\
\leq C(\|q_0\|_{\dot{B}_{2,1}^1})\|\dq\|_{\dot{B}_{2,1}^0},
\end{multline}
and (here we reach the endpoint $d-2=0$ in the remainder)
\begin{multline}
 \|\de G_1\|_{\dot{B}_{2,\infty}^{-1}}=2\|D(u_L+\ub_1).\n\Big(\ln(1+q_1)-\ln(1+q_2)\Big)\|_{\dot{B}_{2,\infty}^{-1}}\\
\leq 2(\|T_D \n\|_{\dot{B}_{2,\infty}^{-1}}+ \|T_\n D\|_{\dot{B}_{2,\infty}^{-1}}+ \|R(D,\n)\|_{\dot{B}_{1,\infty}^{-1}})\\
\leq 2(\|D\|_{L^\infty}\|\n\|_{\dot{B}_{2,\infty}^{-1}} +\|\n\|_{\dot{B}_{\infty,\infty}^{-2}} \|D\|_{\dot{B}_{2,\infty}^1}+ \|D\|_{\dot{B}_{2,2}^1}\|\n\|_{\dot{B}_{2,2}^{-1}})\\
\leq C(\|q_0\|_{\dot{B}_{2,1}^1}) \|u_L+\ub_1\|_{\dot{B}_{2,2}^2 \cap \dot{B}_{\infty,1}^1} \|\dq\|_{\dot{B}_{2,1}^0}.
\end{multline}
Concerning the last term, with the same decomposition, $\de G_2=R_1+R_2$ and when we choose some $\alpha\in]0,1[$ (for $R_2$) we obtain that:
\begin{multline}
 \|R_1\|_{\dot{B}_{2,\infty}^{-1}}\leq \|\dub\|_{\dot{B}_{2,\infty}^1} C(\|q_2\|_{L^\infty})\Big(1+\|q_2\|_{\dot{B}_{2,1}^1}\Big)\|q_2-S_m q_2\|_{\dot{B}_{2,1}^1}\\
\leq C(\|q_0\|_{\dot{B}_{2,1}^1})\Big(\|q_0-S_m q_0\|_{\dot{B}_{2,1}^1} +(e^{CV(t)}-1)\Big)\|\dub\|_{\dot{B}_{2,\infty}^1},
\end{multline}
and
\begin{multline}
 \|R_2\|_{\dot{B}_{2,\infty}^{-1}}\leq \|\dub\|_{\dot{B}_{2,\infty}^{1-\alpha}} \|\ln(1+S_m q_2)\|_{\dot{B}_{2,1}^{1+\alpha}}\leq C(\|q_0\|_{\dot{B}_{2,1}^1}) 2^{\alpha m}\|\dub\|_{\dot{B}_{2,\infty}^{1-\alpha}}.
\end{multline}
As in the previous section, we collect the estimates and obtain:
$$
\|\dq\|_{\Tilde{L}_t^\infty \dot{B}_{2,1}^0} \leq c_0 e^{CV(t)} \|\dub\|_{L_t^1\dot{B}_{2,1}^1},
$$
and with obvious notations,
$$
\beta(t) \leq c_0 e^{CV(t)} \Bigg(\int_0^t (1+V'(\tau))\|\dq\|_{\Tilde{L}_\tau^\infty \dot{B}_{2,1}^0} d\tau +\Big(2^{\alpha m}t^{\frac{\alpha}{2}} +\|q_0-S_m q_0\|_{\dot{B}_{2,1}^1} +(e^{CV(t)}-1)\Big)\beta(t)\Bigg).
$$
under the conditions from the previous section, we get:
$$
\beta(t) \leq c_0 e^{CV(t)} \Bigg(\int_0^t (1+V'(\tau))\|\dub\|_{L_t^1\dot{B}_{2,1}^1} d\tau\Bigg) +\frac{1}{2}\beta(t).
$$
Thanks to Proposition \ref{interpolationlog}, we can use the following logarithmic estimates ($d=2$)
$$
\|\dub\|_{L_t^1\dot{B}_{2,1}^\fd}\leq C \|\dub\|_{L_t^1\dot{B}_{2,\infty}^\fd} \log \Big(e+\frac{\|\dub\|_{L_t^1\dot{B}_{2,\infty}^{\fd-1}} +\|\dub\|_{L_t^1\dot{B}_{2,\infty}^{\fd+1}} }{\|\dub\|_{L_t^1\dot{B}_{2,\infty}^\fd}}\Big)
$$
As $\du=\ub_1-\ub_2$, we can write:
$$
\|\dub\|_{L_t^1\dot{B}_{2,\infty}^0} +\|\dub\|_{L_t^1\dot{B}_{2,\infty}^2}\leq W(t)=W_1(t)+W_2(t),
$$
with $W_i(t)= \|\ub_i\|_{L_t^1\dot{B}_{2,\infty}^0} +\|\ub_i\|_{L_t^1\dot{B}_{2,\infty}^2}$. This function is bounded on $[0,T]$ and the estimates turn into:
$$
\beta(t)\leq C_T \int_0^t (1+V'(\tau)) \beta(\tau) \log \Big(e+\frac{W(T)}{\beta(\tau)}\Big) d\tau.
$$
As we have,
$$
\int_0^1 \frac{dr}{r \log (e+\frac{W(T)}{r})}=\infty,
$$
The Osgood lemma allows us to conclude that $\beta(t)=0$ for all $t\in[0,T]$ (we refer for example to \cite{BCD}, section $3.1.1$). Then the density fluctuation is also zero on this intervall. Then the conclusion is the same as in the case $d\geq 3$.

\section{Global well-posedness}
\label{section5}

In this section we are interested in proving the global well-posedness of (\ref{0.1}) when we assume smallness on the initial data. The proof follows the same lines than in the sections \ref{section3} and \ref{section4}. The only difficulty consists in getting damped effects on the density
in order to deal with the pressure in the remainder terms. To do this we have just to use the estimates in Besov spaces from \cite{DG} or \cite{CD} on the following linear system associated to (\ref{0.1}):
\begin{equation}
\begin{aligned}
&\p_{t}q+v\cdot\n q+{\rm div}u=F,\\
&\p_{t}u+v\cdot\n u-\D u+\n q=G,\\
\end{aligned}
\label{lineaire}
\end{equation}
There, they exhibit the parabolic smoothing effect on $u$ and on
the low frequencies of
$q$, and a damping effect on the high frequencies of $q$. To do this, the authors need to introduce a paralinearisation in order to deal with the convection terms $u\cdot\n q$. More precisely they obtain the following proposition:
\begin{prop}\sl{
\label{5linear1} Let $(q,u)$ a solution of the system (\ref{lineaire}) on
$[0,T[$ , $1-\N<s< 1+\N$ and
$V(t)=\int^{t}_{0}\|\nabla v(\tau)\|_{L^{\infty}\cap \dot{B}^{\N}_{2,2}}d\tau$. We have then the
following estimate for any $T>0$:
$$
\begin{aligned}
&\|(q,u)\|_{\widetilde{L}^{\infty}_{T}(\widetilde{\dot{B}}^{s-1,s}_{2,1})\times \widetilde{L}^{\infty}_{T}(\dot{B}^{s-1}_{2,1})}+\|(q,u)\|_{\widetilde{L}^{1}_{T}(\widetilde{\dot{B}}^{s+1,s}_{2,1})\times \widetilde{L}^{\infty}_{T}(\dot{B}^{s+1}_{2,1})}\\
&\hspace{1cm}\leq Ce^{CV(t)}\big(\|
(q_{0},u_{0})\|_{\widetilde{\dot{B}}^{s-1,s}\times B^{s-1}}+\int^{T}_{0}
e^{-CV(\tau)}\|
(F,G)(\tau)\|_{\widetilde{\dot{B}}^{s-1,s}\times \dot{B}^{s-1}}d\tau\big),\\
\end{aligned}
$$
where $C$ depends only on  $N$ and $s$ and $\widetilde{\dot{B}}^{s_1,s_2}_{2,1}$ denotes the hybrid Besov space with regularity $s_1$ for low frequencies and $s_2$ for high frequencies (we refer to \cite{DG} or \cite{CD} for details).}
\end{prop}
Using this on $\ub$, the rest of the proof simply consists in treating the remainder terms as in section \ref{section3}, for the uniqueness the method follows the same approach as in section \ref{section4}.

\section{Proof of theorem \ref{theo2}}
\label{section6}
Our method from section \ref{section3} and \ref{section4} may be adapted to the study of incompressible density dependent Navier-Stokes equations. This is just a matter of replacing the parabolic model below by a nonstationary Stokes system. More precisely we define $u_{L}$ as the solution of the following system:
\begin{equation}
\begin{aligned}
&\p_{t}u_{L}-\D u_{L}+\n\Pi_{L}=0,\\
&{\rm div}u_{L}=0,\\
&(u_{L})_{\ t=0}=u_{0}.
\end{aligned}
\end{equation}
In the same way than in section \ref{section3}, we are searching solution of the form $u=u_{L}+\bar{u}$ with:
\begin{equation}
\begin{cases}
\begin{aligned}
&\p_t q+(\ub+u_L).\n q +(1+q)\div (\ub+u_L)=0,\\
&\p_t \ub-\Delta \ub+ (\ub+u_L).\n \ub+ \ub.\n u_L+u_L.\n u_L & \\
&\hspace{3cm}-2 D(\ub+u_L).\n\Big(\ln(1+q)\Big)+\frac{1}{1+q}\n (\bar{\Pi}+\Pi_{L})=0,\\
&{\rm div}\bar{u}=0.
\end{aligned}
\end{cases}
\label{systbarre1}
\end{equation}
By applying the operator ${\rm curl}$ to the momentum equation, we obtain that:
\begin{equation}
\begin{cases}
\begin{aligned}
&\p_t q+(\ub+u_L).\n q +(1+q)\div (\ub+u_L)=0,\\
&\p_t{\rm curl} \ub-\Delta {\rm curl}\ub+ {\rm curl}\biggl(\ub+u_L).\n \ub+ \ub.\n u_L+u_L.\n u_L & \\
&\hspace{2cm}-2 D(\ub+u_L).\n\Big(\ln(1+q)\Big)\biggl)+\nabla(\frac{1}{1+q}):\n (\bar{\Pi}+\Pi_{L})=0,\\
&{\rm div}\bar{u}=0.
\end{aligned}
\end{cases}
\label{systbarre1}
\end{equation}
By following the same idea as in section \ref{section3}, we are able to estimate $\bar{u}$ in $\widetilde{L}^{\infty}_{T}(\dot{B}^{\N-1}_{2,1})\cap \widetilde{L}^{1}_{T}(\dot{B}^{\N+1}_{2,1})$ by proving estimate on ${\rm curl}u$  in  $\widetilde{L}^{\infty}_{T}(\dot{B}^{\N-2}_{2,1})\cap \widetilde{L}^{1}_{T}(\dot{B}^{\N}_{2,1})$, in order to deal with the pressure $\bar{\Pi}$ it is sufficent to adapt the idea of \cite{AP,H} as $\bar{\Pi}$ verifies an elliptic equation.
 \\
\\
\textbf{Aknowledgements:} The authors would like to thank Rapha\"el Danchin and Pierre Germain for many fruitful discussions.
 
\end{document}